\documentclass[12pt]{article}

\usepackage{amssymb,amsmath,amsfonts}
\usepackage[T2A]{fontenc}
\usepackage[cp1251]{inputenc}

\begin{document}

%%%%%%%%%%%%%%%%%%%%%%%%%%%%%%%%%%%%%%%%%%%%%%%%%%%%%%%%%%%%%%%%%%%%%%%%%%%%%%

\begin{center}
\textbf{Symmetries of the Black-Scholes-Merton equation for European options}\\
L.N.~Bakirova$^1$, M.A.~Shurygina$^2$, V.V.~Shurygin, jr.$^3$
\end{center}

\begin{flushleft}
$^1$~ICL Group, Kazan, Russia, Lbakir@mail.ru\\
$^2$~Kazan National Research Technical University named after A.N. Tupolev--KAI, Kazan, Russia, m-akm@mail.ru\\
$^3$~N.I. Lobachevskii Institute of Mathematics and Mechanics, Kazan Federal University, Kazan, Russia, vshjr@yandex.ru, corresponding author
\end{flushleft}

\begin{abstract}
The aim of the present paper is the clarification of the result of  
A.\,Paliathanasis, K.\,Krishnakumar, K.M.\,Tamizhmani
and
P.G.L.\,Leach on the symmetry Lie algebra of the Black-Scholes-Merton
equation for European options.
\end{abstract}



\section{Main result}
%%%%%%%%%%%%%%%%%%%%%%%%%%%%%%%%%%%%%%%

The  Black-Scholes-Merton (BSM)~\cite{BS, BS1, M} model is one of the most important concepts in modern financial theory. 
It is used for the valuation of stock options, taking into account the impact of time and other risk factors. 

The classical  BSM  model is described by a second order PDE
$$
u_t+\frac12\sigma^2x^2u_{xx}+rxu_x-ru=0.
$$
In recent years the number of papers is devoted to the determining the Lie algebra of symmetries of
some PDEs that generalize this model.
First of all we mention the papers of
R.K.\,Gazizov, and N.H.\,Ibragimov~\cite{GI} and
O.\,Sinkala,  P.G.L.\,Leach, J.\,O’Hara~\cite{SLH}. 
The paper of Y.\,Bozhkov and S.\,Dimas~\cite{BD}
solves the problem of group classification of the generalized BSM equation$$
u_t+\frac12\sigma^2x^2u_{xx}+rxu_x+f(u)=0,
$$
where $f(u)$ is an arbitrary smooth function.

A.\,Paliathanasis, K.\,Krishnakumar, K.M.\,Tamizhmani, P.G.L.\,Leach~\cite{PKTL}
considered the BSM equation for  European options
with stochastic volatility for which the premium term depends only upon the return-to-risk ratio.
It has the form
\begin{equation}
\label{BSMe1}
 \dfrac{1}{2}f^2(y)x^2u_{xx}+\rho\beta xf(y)u_{xy}+\dfrac{1}{2}\beta^2u_{yy}
+rxu_{x}+\left(\alpha(m-y)-\beta\rho\dfrac{\mu-r}{f(y)}\right)u_{y}-ru+u_t=0,
\end{equation}
where
$t$, $x$, $y$ are independent variables,  $u = u (t, x, y)$  is the value of the option, $f(y)$ is an arbitrary smooth function, and $r$, $\rho$, $m$, $\mu$, $\alpha$, $\beta$ are real parameters with $|\rho|<1$.

The authors calculate the Lie algebra of symmetries of this equation in the case when $f(y)={\rm const}$, 
and claim that in the case of a non-constant function
$f(y)$ this Lie algebra is the direct sum of a three-dimensional commutative Lie algebra and the infinite dimensional
commutative Lie algebra which corresponds to the solutions on the equation.
It should be mentioned that there are some misprints in their  paper.

In the present paper we  add some corrections to this result.
We show that this PDE admits additional symmetries also in the case when
$$
f(y)=\frac{k}{y-m}, ~~{k}={\rm const}.
$$
For this case we denote 
$$
g=2\left(\alpha+\dfrac{\rho \beta(\mu-r)}{k}\right).
$$
The computations  were performed
using the {\tt Maple} packages {\tt DifferentialGeometry} and {\tt JetCalculus}
by I.M.\,Anderson.

\medskip
\textbf{Theorem.}
{\it For arbitrary function $f(y)$ the equation {\rm(\ref{BSMe1})} admits the Lie symmetries
$$
X_1=\dfrac{\partial}{\partial t},~~ X_2=x\dfrac{\partial}{\partial
x},~~ X_3=u\dfrac{\partial}{\partial u},~~ X_b=b(x,\, y,\,
t)\dfrac{\partial}{\partial u},
$$
where $b(x,y,t)$ is the solution of {\rm(\ref{BSMe1})}. 

Moreover, in the case when
$f=f_0={\rm const}$, it admits three additional symmetries
\begin{gather*}
 X_4=e^{-\alpha t}\dfrac{\partial}{\partial y},
%\end{equation*}
%\begin{equation*}
\\
X_5=f_0^2\left(\rho^2+\alpha t\right)x\dfrac{\partial}{\partial x}+f_0\rho \beta \dfrac{\partial}{\partial y}+
\dfrac{1}{2}\alpha \left(t (f_0^2-2r)+2 \ln x\right)u\dfrac{\partial}{\partial u},
%\end{equation*}
%\begin{equation*}
\\
X_6=2 \beta f_0^2\rho e^{\alpha t} x\dfrac{\partial}{\partial x}+\beta^2 f_0 e^{\alpha t} \dfrac{\partial}{\partial y}-
2 e^{\alpha t}\left(\alpha f_0 (m-y) + \beta \rho(r-\mu)\right)u \dfrac{\partial}{\partial u}.
\end{gather*}

In the case 
$f(y)=\dfrac{k}{y-m}$ and $g\ne0$ 
the equation {\rm(\ref{BSMe1})}  takes the form
$$
 \dfrac{1}{2}f^2(y)x^2u_{xx}+\rho\beta xf(y)u_{xy}+\dfrac{1}{2}\beta^2u_{yy}
+rxu_{x}+\frac12g(m-y)u_{y}-ru+u_t=0
$$
and admits two additional symmetries
\begin{equation*}
%X_4=e^{-gt} r x \dfrac{\partial}{\partial x}+\dfrac{1}{2} g e^{-gt}(m-y) \dfrac{\partial}{\partial y}+e^{-gt} \dfrac{\partial}{\partial t}+u r e^{-gt} \dfrac{\partial}{\partial u},
X_4=e^{-gt}\left( r x \dfrac{\partial}{\partial x}-\dfrac{1}{2} g (y-m) \dfrac{\partial}{\partial y}+ \dfrac{\partial}{\partial t}+u r  \dfrac{\partial}{\partial u}\right),
\end{equation*}
\begin{equation*}
X_5=e^{gt}\left(\dfrac{x }{\beta}\left(\rho g k +\beta r\right)\dfrac{\partial}{\partial x}+\dfrac{1}{2}g (y-m)\dfrac{\partial}{\partial y}
+\dfrac{\partial}{\partial t}+\dfrac{u }{2\beta^2}
\left(g^2 (y-m)^2+\beta^2 (2r-g)\right)\dfrac{\partial}{\partial u}\right).
\end{equation*}
In the case 
$f(y)=\dfrac{k}{y-m}$ and $g=0$ 
the equation {\rm(\ref{BSMe1})} takes the form
$$
 \dfrac{1}{2}f^2(y)x^2u_{xx}+\rho\beta xf(y)u_{xy}+\dfrac{1}{2}\beta^2u_{yy}
+rxu_{x}-ru+u_t=0
$$
and admits two additional symmetries
\begin{equation*}
X_4=\dfrac{1}{2\beta}\left(\rho  k +2 \beta rt \right)x
\dfrac{\partial}{\partial
x}+\dfrac{y-m}{2}\dfrac{\partial}{\partial y}+
t\dfrac{\partial}{\partial t}+rtu \dfrac{\partial}{\partial u},
\end{equation*}
\begin{equation*}
X_5=\dfrac{xt}{\beta}\left(\rho k + \beta rt \right) \dfrac{\partial}{\partial x}
+t\left(y-m\right)\dfrac{\partial}{\partial y}
+t^2 \dfrac{\partial}{\partial t}+\dfrac{u}{2\beta^2}\left((y-m)^2 +
\beta^2 (2r t^2-t)\right)\dfrac{\partial}{\partial u}.
\end{equation*}
}

\section{Invariant solutions}

In this section, we apply the Lie symmetries in order to reduce the equation (\ref{BSMe1}) and to construct the invariant solutions.
To do this, we need to use 2-dimensional Lie subalgebras.
As it is mentioned in~\cite{PKTL}, solutions in which $u$ does not depend upon
some of the independent variables, are not interesting.
Therefore, the reductions are performed with modified symmetry vectors like $ X_1 + k_1X_3$, 
$ X_2 + k_2X_3$ and some other.

\subsection{The case $f=f_0={\rm const}$}

In this subsection we repeat the results of A.\,Paliathanasis et al. with minor corrections.
The Lie Brackets of the Lie algebra are given in the table.

\bigskip

\begin{tabular}{ccccccc}
\hline
$[X_i,X_j]$ & $X_1$ & $X_2$ & $X_3$ & $X_4$ & $X_5$ & $X_6$ \\
\hline
$X_1$ & 0 & 0 & 0 & $-\alpha X_4$  & $ f_0^2\alpha X_2+\alpha(\frac12f_0^2-r)X_3$ & $\alpha X_6$ \\
$X_2$ & 0 & 0 & 0 & 0 & $\alpha X_3$ & 0 \\
$X_3$ & 0 & 0 & 0 & 0 & 0 & 0 \\
$X_4$ & $\alpha X_4$ & 0 & 0 & 0 & 0 & $2\alpha f_0 X_3$ \\
$X_5$ & $ -f_0^2\alpha X_2-\alpha(\frac12f_0^2-r)X_3$ & $-\alpha X_3$ & 0 & 0 & 0 & 0 \\
$X_6$ & $-\alpha X_6$ & 0 & 0 & $-2\alpha f_0 X_3$ & 0 & 0 \\
\hline
\end{tabular}

\bigskip

Throughout this subsection we denote
$Y_1 = X_1 + k_1X_3$, $Y_2 = X_2 + k_2X_3$.

\smallskip
1) Algebra $\{Y_1,Y_2\}$.

The invariant solution  has the form
$u(x,y,t) = e^{k_1t}x^{k_2}w(y) $, where
$w(y)$ satisfies the equation
$$
\beta^2f_0w''+
2(\rho\beta f_0^2k_2+\alpha(m-y)f_0+\beta\rho(r-\mu))w'+
((k_2-1)(f_0^2k_2+2r)+2k_1)w=0.
$$
Its solution is expressed via Kummer functions $M$, $U$, so finally we get:
$$
u(x,y,t) = e^{k_1t}x^{k_2}(C_1 M(\gamma,\tfrac12,z) +C_2U(\gamma,\tfrac12,z)),
$$
where
$$
z=\frac{(k_2\rho\beta f_0^2+\beta\rho(r-\mu)-\alpha f_0(y-m))^2}{\alpha\beta^2 f_0^2},
~~
\gamma=\frac{k_2(1-k_2)f_0^2+2r(1-k_2)-2k_1}{4\alpha}.
$$

2) Algebra $\{Y_2, X_4+kX_3\}$.

The invariant solution  has the form
$u(x,y,t)=\exp(ke^{\alpha t}y)x^{k_2}w(t)$, where
$w(t)$ satisfies the equation
$$
2f_0w' +
%(k_2(k_2-1)f_0^3 + 2kk_2e^{\alpha t}\rho\beta f_0^2
%+(\beta^2k^2e^{2\alpha t}+2ke^{\alpha t}\alpha m+2r(k_2-1))f_0 
%-2e^{\alpha t}k\rho\beta(\mu-r) )
\left(\beta^2k^2f_0e^{2\alpha t}  +  
2k\beta\rho\left(k_2 f_0^2+\frac{\alpha m f_0}{\rho\beta} -(\mu-r) \right)e^{\alpha t}
+ (k_2-1)f_0( k_2f_0^2  +2r) \right)
w=0.
$$
Its solution is
$$
w(t)=C\exp\left(-\frac{\beta^2k^2}{4\alpha} e^{2\alpha t} 
-\frac{k\beta\rho}{\alpha f_0} e^{\alpha t} \left( f_0^2k_2+\frac{\alpha f_0 m}{\beta\rho}-\mu+r\right)
-\frac12 t(k_2-1)(f_0^2k_2+2r)\right).
$$

3) Algebra $\{Y_2, X_6+kX_3\}$.

The invariant solution  has the form
$$
u(x,y,t)=\exp\left(\dfrac{\alpha}{\beta^2}y^2
-\left(\dfrac{2\alpha f_0m+ke^{-\alpha t}}{\beta^2 f_0}
 + \dfrac{2\rho(r-\mu+f_0^2k_2)}{\beta f_0}\right)y  \right)  x^{k_2} w(t),
$$ where
$w(t)$ satisfies the equation
$$
2\beta^2 f_0^2 w' +\left( f_0^2\beta^2(2\alpha+(k_2-1)(f_0^2k_2+2r))
+ 2k(k_2\beta f_0^2\rho+\alpha f_0m +\beta\rho(r-\mu)) e^{-\alpha t}
+k^2e^{-2\alpha t}\right)w=0.
$$
Its solution is
\begin{multline*}
w(t)=C\exp\Bigl(
\frac1{4\alpha\beta^2f_0^2}
\bigl( 
k^2e^{-2\alpha t}
+4k(k_2\beta f_0^2\rho+\alpha f_0m+\beta\rho(r-\mu))e^{-\alpha t}-\\
-2t\alpha\beta^2f_0^2(2\alpha+(k_2-1)(f_0^2k_2+2r))
\bigr)
\Bigr).
\end{multline*}

4) Algebra $\{X_5, X_4+kX_3\}$.

The invariant solution  has the form
$$u(x,y,t)=
\exp\left(kye^{\alpha t} +\frac{\alpha \ln^2x}{2f_0^2(\rho^2+\alpha t)}\right) x^{\psi(t)}
w(t),
$$
where
$$
\psi(t)=
-\frac{2k\beta\rho f_0e^{\alpha t}+(2r-f_0^2)t\alpha}{2f_0^2(\rho^2+\alpha t)}.
$$
The function $w(t)$ is the solution of the linear first order ODE
\begin{multline*}
8f_0^2(\rho^2+\alpha t)^2 w'
+
\Bigl(4e^{2\alpha t}k^2f_0^2\beta^2(\alpha^2t^2+\rho^2-\rho^4)\,+
\\
+4e^{\alpha t}\bigl(
kf_0\alpha^2(\beta f_0^2\rho+2f_0\alpha m-2\rho\beta\mu)t^2+
k\alpha\rho^2f_0(\beta f_0^2\rho+4\alpha f_0 m+2\beta\rho(r-2\mu))t
\,+\\
+(2k f_0^2\alpha m\rho^4  -2k\beta\rho^3f_0(\rho^2(\mu-r)+r) -1 )
\bigr)
-
\\
-
\bigl(\alpha^2(2r+f_0^2)^2t^2+
2\alpha(f_0^4\rho^2+2f_0^2(2r\rho^2-\alpha)+4r^2\rho^2)t
+4f_0^2\rho^2(2r\rho^2-\alpha)
\bigr)
\Bigr)w=0.
\end{multline*}
The final form of $w$ is too cumbersome to give it here.

5) Algebra $\{X_5, X_6\}$.

The invariant solution  has the form
$$u(x,y,t)=
 x^{\psi(t)} \exp(\varphi(x,y,t))
w(t),
$$
where we denoted
$$
\psi(t)=\frac{\bigl(t(r-\frac12f_0^2)\beta-2\rho f_0(m-y)\bigr)\alpha+2\beta\rho^2(\mu-r)}{f_0^2(\rho^2-\alpha t)\beta},
$$
$$
\varphi(x,y,t)=\alpha\beta^2\ln^2x
-2yf_0\bigl(
(\rho^2+\alpha t)(2m-y)\alpha f_0
+\alpha\beta\rho(f_0^2-2\mu)t + 2\rho^3\beta(r-\mu)
\bigr).
$$
The function $w(t)$ has the form
$$
w(t)=
\frac{C\exp\left(\dfrac{\xi(t)}{\alpha \beta^2 f_0^2(\alpha t-\rho^2)} \right)}{\sqrt{\alpha t-\rho^2}},
$$
where
$$
\xi(t)=\frac18 \Bigl(\alpha\beta^2(\alpha t-\rho^2)\bigl((2r+f_0^2)^2-8f_0^2\alpha\bigr)
+\rho^2\bigl(\beta\rho (f_0^2-4\mu+2r)+ 4\alpha f_0m\bigr)^2\Bigr).
$$

%%%%%%%%%%%%%

\subsection{The case $f=\dfrac{k}{y-m}$, $g\ne0$}

For convenience we will use $g=2\left(\alpha+\dfrac{\rho \beta(\mu-r)}{k}\right)$ instead of $\alpha$.
The Lie Brackets of the Lie algebra are given in the table.

\bigskip

\begin{tabular}{cccccc}
\hline
$[X_i,X_j]$ & $X_1$ & $X_2$ & $X_3$ & $X_4$ & $X_5$  \\
\hline
$X_1$ & 0 & 0 & 0 & $-gX_4$  & $gX_5$  \\
$X_2$ & 0 & 0 & 0 & 0  & 0 \\
$X_3$ & 0 & 0 & 0 & 0 & 0  \\
$X_4$ & $gX_4$ & 0 & 0 & 0 &  $X$ \\
$X_5$ & $-gX_5$  & 0 & 0 & $-X$ & 0 \\
\hline
\end{tabular}

\bigskip
where we denoted
$$
X=2gX_1+\Bigl(\frac{g^2\rho k}{\beta}+2 r\Bigr) X_2-\frac12g(g-4r)X_3.
$$

1) Algebra $\{X_1 + aX_3, X_2 + bX_3\}$.

The invariant solution  has the form
$$
u(x,y,t) = e^{at}x^bw(y),
$$
where $w(y)$ satisfies the  second order ODE
$$
\beta^2(y-m)^2 w'' 
+(2\beta\rho bk -g(y-m)^2)(y-m)w'
+\bigl(2(a+r(b-1))(y-m)^2+k^2b(b-1)\bigr)w=0.
$$
Its solution can be expressed via Whittaker functions $M_{\kappa\nu}$, $W_{\kappa\nu}$:
$$
w(y)=\exp\left(\frac{g(y-m)^2}{4\beta^2}\right)(y-m)^{\gamma}
\left(
M_{\kappa\nu}\left(\frac{g(y-m)^2}{2\beta^2}\right)   + W_{\kappa\nu}\left(\frac{g(y-m)^2}{2\beta^2}\right)
\right),
~~\gamma=-\frac12+\frac{kb\rho}\beta,
$$
and
$$
\kappa=\frac{(\beta+2kb\rho)g+4\beta(a+r(b-1))}{4g\beta} ,
~~
\nu=\frac1{4\beta}\sqrt{\beta^2 +4k(k-\rho\beta)b -4k^2b^2 (1-\rho^2)}.
$$

2) Algebra $\{ X_1 + aX_3, X_4\}$

The invariant solution  has the form
$$
u(x,y,t) = e^{at}x^{1-a/r}w(h(t,y)),
~~
h(t,y)=x^g(y-m)^{2r}.
$$
The function $w(h)$ satisfies the second order Euler equation
\begin{multline*}
r^2(4\beta^2 r^2+g^2k^2+4rkg\rho\beta)h^2w''-\\
-r(k^2g(2a-g(r+1)) +4r\rho\beta(a-g(r+1))k  -2r^2\beta^2(2r-1))hw'+ak^2(a-r)w=0.
\end{multline*}

3) Algebra $\{ X_2 + bX_3, X_4\}$

The invariant solution  has the form
$$
u(x,y,t) =x^b(y-m)^{\gamma}
w(h(t,y)),
~~
h(t,y)=e^{gt}(y-m)^{2},
~~
\gamma=\frac{2r(b-1)}g.
$$
The function $w(h)$ satisfies the second order Euler equation
\begin{multline*}
4\beta^2 g^2 h^2w''
+2g\beta(g(2bk\rho+\beta)+4r\beta(b-1)) hw'+\\
+(b(k^2g^2+4r^2\beta^2+4kgr\rho\beta)-2r\beta^2(2r+g))w=0.
\end{multline*}

4) Algebra $\{X_2+bX_3, X_5+pX_3\}$.

The invariant solution  has the form
$$
u(x,y,t) = x^b(y-m)^{\gamma}\exp(\varphi(t,y))w(h(t,y)),
$$
where we denoted
$$ 
h(t,y)=(y-m)^2e^{-gt}, ~~
\varphi(t,y)=-\frac{g^2y(2m-y)+2p\beta^2e^{-gt}}{2g\beta^2};~~
 \gamma=-\frac{(g-2r+2rb)\beta+2b\rho kg}{g\beta}.
$$

The function $w(h)$ satisfies the second order ODE
\begin{multline*}
4g^2\beta^2h^2w''-2\beta g(g\beta+2b\rho kg+4\beta r(b-1))hw'+ (2pg^2h +\\
+6g\beta^2 rb-8\beta^2 r^2b+4\beta^2 r^2b^2-6g\beta^2 r+k^2b^2g^2-k^2bg^2+4b^2\rho \beta kg r+\\
+2g^2\beta^2+4\beta^2 r^2+4g^2\beta b\rho k-4b\rho\beta rkg)w=0.
\end{multline*}
Its solution is
$$
w(h)=h^\gamma\left(C_1J_n \left(\frac{\sqrt{2ph}}\beta \right) +C_2Y_n \left(\frac{\sqrt{2ph}}\beta \right)\right),
$$
where
$$\gamma=\frac{(4r(b-1)+3g)\beta+2bkg\rho}{4g\beta},
~~n=\frac1{2\beta}\sqrt{(4(\rho^2-1)k^2b^2-4k(4\beta\rho-k)b+\beta^2}
$$
and $J_n$, $Y_n$ are the Bessel functions.

5) Algebra $\{X_2+aX_1, X_4\}$.

The invariant solution  has the form
$$
u(x,y,t)=(y-m)^{-2r/g}w(h(x,y,t)),
~~
h(x,y,t)=e^{gt}x^{-ag}(y-m)^{2(1-ar)}.
$$
The function $w(h)$ satisfies the second order Euler equation
\begin{multline*}
g^2((4\beta^2 r^2+4kgr\rho\beta+k^2g^2)a^2   -4\beta(gk\rho+2\beta r)a  +4\beta^2)h^2w''+
\\
+g(g(4\beta^2r^2+4kgr\rho\beta+k^2g^2)a^2 + (k^2g^2+8\beta^2r^2+4kgr\rho\beta-6gr\beta^2-4kg^2\rho\beta)a +2g\beta^2)  hw'
+\\
+2r\beta^2(2r+g)w=0.
\end{multline*}

6) Algebra $\{X_2+aX_1, X_5\}$.

The invariant solution  has the form
$$
u(x,y,t)=(y-m)^{-1+2r/g}\exp\left(\frac{gy(y-2m)}{2\beta^2}\right)w(h(x,y,t)),
$$
where
$$
h(x,y,t)=\frac{(y-m)^\gamma e^{\beta gt}}{x^{a\beta g}},
~~
\gamma=2(agk\rho+ar\beta -\beta).
$$
The function $w(h)$ satisfies the second order Euler equation
\begin{multline*}
g^4\beta^2(a^2(4r^2\beta^2+4rgk\beta\rho+k^2g^2) -4a\beta(gk\rho+2\beta r)+4\beta^2)h^2w''+
\\
+g^3\beta(4g(ar-1)^2\beta^3 +2(ar-1)(2ag^2k\rho+4r-3g)\beta^2 +agk(akg^2+4\rho(r-g))\beta +k^2g^2a )  hw'
+\\
+(g-r)(g-2r)w=0.
\end{multline*}
%%%%%%%%%%%%

\subsection{The case $f=\dfrac{k}{y-m}$, $g=0$}

The Lie Brackets of the Lie algebra are given in the table.

\bigskip

\begin{tabular}{cccccc}
\hline
$[X_i,X_j]$ & $X_1$ & $X_2$ & $X_3$ & $X_4$ & $X_5$  \\
\hline
$X_1$ & 0 & 0 & 0 & $X_1+rX_2+rX_3$  & $-\frac12X_3+2X_4$  \\
$X_2$ & 0 & 0 & 0 & 0  & 0 \\
$X_3$ & 0 & 0 & 0 & 0 & 0  \\
$X_4$ & $-X_1-rX_2-rX_3$ & 0 & 0 & 0 &  $X_5$ \\
$X_5$ & $\frac12X_3-2X_4$  & 0 & 0 & $-X_5$ & 0 \\
\hline
\end{tabular}

\bigskip
%For convenience we will use $g=2\left(\alpha+\dfrac{\rho \beta(\mu-r)}{k}\right)$ instead of $\alpha$.

1) Algebra $\{X_1 + aX_3, X_2 + bX_3\}$.

The invariant solution  has the form
$$
u(x,y,t) = e^{at}x^bw(y),
$$
where $w(y)$ satisfies the  second order ODE
$$
\beta^2(y-m)^2 w'' 
+2\beta\rho bk(y-m) w'
+\bigl(2(a+r(b-1))(y-m)^2+k^2b(b-1)\bigr)w=0.
$$
Its solution can be expressed via Whittaker functions:
$$
w(y)=(y-m)^{-bk\rho/\beta}(C_1M_{0,\nu}(z)+C_2W_{0,\nu}(z)), 
$$
where
$$
z=\frac{2\sqrt{(2-2b)r-2a}(y-m)}{\beta},
~~
\nu=\frac{\sqrt{4b(1-b+b\rho^2)k^2+\beta(\beta-4bk\rho)}}{2\beta}.
$$

2) Algebra $\{ X_2 + bX_3, X_4\}$.

The invariant solution  has the form
$$
u(x,y,t) = x^b  (y-m)^{\gamma}\exp(\varphi(t,y))
w(h(t,y)),
$$
where we denoted
$$ 
h(t,y)=\frac t{(y-m)^2},
~~
\varphi(t,y)=\frac{r(b-1)(2m-y)ty}{(y-m)^2},
~~ \gamma=-\frac{b\rho k}{\beta}.
$$

The function $w(h)$ satisfies the second order ODE
\begin{multline*}
4\beta^2h^2w'' + 2(4rm^2\beta^2(b-1)h^2+3\beta^2 h+1)w'+\\
+
\bigl(
4m^4r^2\beta^2(b-1)^2 h^2 + 6rm^2\beta^2(b-1)h
+(k^2b+2m^2r)(b-1)+bk\rho (\beta-bk\rho)
\bigr)w=0.
\end{multline*}
Its solution can be expressed via Kummer  functions:
$$
w(h)=e^{-m^2hr(b-1)}h^{-\gamma}\left( C_1 M\left(\gamma, 2\gamma+\frac12, \frac1{2\beta^2h}\right)
+C_1 U\left(\gamma, 2\gamma+\frac12, \frac1{2\beta^2h}\right)
\right),
$$
where
$$
\gamma=\frac {\beta+\sqrt{4(\rho^2-1)k^2b^2+4k(k-4\rho\beta)b+\beta^2}}{4\beta}.
$$

3) Algebra $\{ X_2 + bX_3, X_5\}$.

The invariant solution  has the form
$$
u(x,y,t) = x^b  (y-m)^{\gamma}\exp(\varphi(t,y))
w(h(t,y)),
$$
where we denoted
$$ 
h(t,y)=\frac t{y-m},
~~
\varphi(t,y)=\frac{y(y-m)}{2\beta t} -\frac{r(b-1)ty}{y-m},
~~ \gamma=-\frac12-\frac{b\rho k}{\beta}.
$$
The function $w(h)$ satisfies the second order ODE
\begin{multline*}
4\beta^4h^4w''
-4\beta^2h^2(2\beta^2 rm(b-1)h^2 -3\beta^2 h+ m)w'
+\bigl(4\beta^4 r^2m^2(b-1)^2 h^4
-12\beta^4rm(b-1)h^3+
\\
+\beta^2( 4kb(k(b-1)-kb\rho^2-\rho\beta)+ 4m^2r(b-1)+3\beta^2)h^2
-2\beta^2mh+m^2
\bigr)w=0.
\end{multline*}
Its solution is
$$
w(h)=(C_1h^{\lambda_1}+C_2h^{\lambda_2})\exp\left(\frac{m(2r\beta^2(b-1)h^2-1)}{\beta^2 h} \right),
$$
where $\lambda_1$, $\lambda_2$ are the roots of the quadratic equation
$$
\lambda^2+ 2\lambda
+\frac1{4\beta^2}\bigl(4(1-\rho^2)b^2k^2+3\beta^2-4k^2b+4bk\beta\rho\bigr)=0.
$$

4) Algebra $\{ X_4, X_5\}$.

The invariant solution  has the form
$$
u=x\sqrt{rt}\exp\left(\frac{(y-m)^2}{2\beta^2t}\right)(y-m)^{-1-k\rho/\beta}
w(h),
$$
where
%$\rho=\sin\theta$, 
$h(t,x,y)=k\rho\ln(y-m)+\beta(rt-\ln x)$, and the function $w(h)$ satisfies the second order ODE with the constant coefficients
%$$
%k^2\beta^2(1-\rho^2)\lambda^2-
%k\beta(k+\beta\rho-2k\rho^2)\lambda
%+(2\beta^2 +k\beta\rho-k^2\rho^2)=0.
%$$
$$
k^2\beta^2(1-\rho^2)w''+k\beta(k+\beta\rho-2k\rho^2)w'+(\beta+k\rho)(2\beta-k\rho)w=0.
$$

{\bf Author Contributions:} Vadim V. Shurygin, jr. determined the problem and the methods
for the solution, checked the results and wrote the paper. Landysh N. Bakirova calculated the Lie algebra of symmetries. Marina A. Shurygina computed the invariant solutions. 

{\bf Conflicts of Interest:} The authors declare no conflict of interest.

%%%%%%%%%%%%%%%%%%%%%%%%%%%%%%%%%%%%%%%%%%%%%%%%%%%%%%%%%%%%%%%%%%%%%%%%%%%%%%%%%%%%%%%%%%%

\end{document}